\documentclass[11pt,reqno]{amsart}
\usepackage{geometry} 
\usepackage{xcolor}
\usepackage{soul}
\usepackage{tikz}
\usepackage{tikz-cd}
\usepackage{color}
\usepackage{hyperref}
\hypersetup{
  colorlinks   = true, 
  urlcolor     = blue, 
  linkcolor    = blue, 
  citecolor   = red 
}

\newtheorem{theorem}{Theorem}[section]
\newtheorem{definition}[theorem]{Definition}
\newtheorem{remark}[theorem]{Remark}
\newtheorem{lemma}[theorem]{Lemma}
\newtheorem{proposition}[theorem]{Proposition}

\everymath{\displaystyle}

\newcommand{\leve}{\left\vert}
\newcommand{\rive}{\right\vert}
\newcommand{\R}{\mathrm{R}}        
\newcommand{\Ric}{\mathrm{Ric}}   
\newcommand{\us}{\underline{s}}   
\newcommand{\uC}{\underline{C}}  
\newcommand{\slashg}{g \mkern-8.7mu \slash}  
\newcommand{\tr}{\mathrm{tr}}              
\newcommand{\slashGamma}{\Gamma \mkern-11mu \slash\,}     
\newcommand{\uL}{\underline{L}}

\newcommand{\ueta}{\underline{\eta}}  
\newcommand{\hatchi}{\hat{\chi}}          
\newcommand{\uchi}{\underline{\chi}}   
\newcommand{\hatuchi}{\hat{\uchi}}      
\newcommand{\uomega}{\underline{\omega}}   

\newcommand{\ed}{\mathrm{d}}        
\newcommand{\slashd}{\ed\mkern-8.5mu \slash \,}      
\newcommand{\slashnabla}{\nabla\mkern-12mu \slash \,}
\newcommand{\slashDelta}{\Delta \mkern-12mu \slash \mkern+2mu}        

\newcommand{\sym}{\mathrm{sym}}   

\newcommand{\uf}{\underline{f}}    
\newcommand{\vareps}{\varepsilon}     
\newcommand{\uvarepsilon}{\underline{\varepsilon}}   
\newcommand{\uvareps}{\underline{\varepsilon}}    

\newcommand{\calH}{\mathcal{H}}

\newcommand{\dL}{\dot{L}}
\newcommand{\duL}{\dot{\uL}}
\newcommand{\dpartial}{\dot{\partial}}
\newcommand{\dslashg}{\dot{\slashg}}

\newcommand{\dtr}{\dot{\tr}}
\newcommand{\dchi}{\dot{\chi}}
\newcommand{\duchi}{\dot{\uchi}}
\newcommand{\deta}{\dot{\eta}}

\newcommand{\dslashd}{\dot{\slashd}}

\newcommand{\bslashg}{\bar{\slashg}}
\newcommand{\bchi}{\bar{\chi}}
\newcommand{\btr}{\bar{\tr}}

\title[Marginal Tubes and Foliations by Marginal Surfaces]{Marginal Tubes and Foliations by Marginal Surfaces}
\author{Pengyu Le}
\address{Yanqi Lake Beijing Institute of Mathematical Sciences and Applications, Beijing, China}
\email{pengyu.le@bimsa.cn}
\date{} 

\begin{document}

\begin{abstract}
In this paper, we introduce the notion of a marginal tube, which is a hypersurface foliated by marginal surfaces. It generalises the notion of a marginally trapped tube and several notions of black hole horizons, for example trapping horizons, isolated horizons, dynamical horizons, etc. We prove that if every spacelike section of a marginal tube is a marginal surface, then the marginal tube is null. There is no assumption on the topology of the marginal tube. To prove it, we study the geometry of spacelike surfaces in a $4$-dimensional spacetime with the help of double null coordinate systems. The result is valid for arbitrary 4-dimensional spacetimes, regardless of any energy condition.

\end{abstract}
\maketitle
\tableofcontents

\section{Introduction}\label{section 1}
\noindent
A trapped surface is a spacelike surface where the area element decreases pointwise for every infinitesimal displacement along the future null normal direction. This concept is first introduced by Penrose in \cite{P2}.
\begin{definition}
$\Sigma$ is called trapped if two future null expansions are both negative on $\Sigma$, i.e.
\begin{align} \nonumber \tr \chi_{\Sigma} < 0  \quad   \text{and} \quad \tr \uchi_{\Sigma} < 0. \end{align}
\end{definition}

Based on this concept, Penrose proved his famous incompleteness theorem, that a spacetime, with a non-compact Cauchy hypersurface and satisfying the null convergence condition that $\Ric(X,X) \geq 0$ for all null vectors $X$, is future null geodesically incomplete if it contains a closed trapped surface. Then with the notion of null infinity introduced in \cite{P1}, Penrose showed that if the above mentioned spacetime has a complete future null infinity, then it must contain a black hole whose interior contains the closed trapped surface.

The result of Penrose is generalised in \cite{HE} to partial Cauchy hypersurfaces and outer trapped surfaces. For a closed surface in a partial Cauchy hypersurface, one can distinguish its future outgoing and incoming null normal directions. An outer trapped surface is an orientable closed spacelike surface where the future outgoing null normal geodesics are converging, where it does not matter the incoming null normal geodesics are converging or not.

Penrose's result and its generalisation in \cite{HE} show that both closed trapped surfaces and outer trapped surfaces are in the black hole. It is convenient to define the boundary of the black hole to be the future boundary of the past of the future null infinity, which is called the event horizon. The event horizon is a global concept, which can not be located without the whole future behaviour of the spacetime.

\cite{HE} introduces the trapped region as the set of all points which an outer trapped surface passes through. The outer boundary of a connected component of the trapped region is called an apparent horizon. It is a 2-dimensional surface. The apparent horizon is an appropriate sort of horizon which depends only on the local geometry of the spacetime on the surface. It is showed in \cite{HE} that an apparent horizon is marginally outer trapped, i.e. the future outgoing null expansion is zero.

The definition of a apparent horizon in \cite{HE} inspires \cite{Ha} to introduce the trapping boundary. In a spacetime, the inside and outside of a 2-dimensional surface cannot be distinguished in general, therefore \cite{Ha} defines a trapped region of a spacetime to be a connected subset of the set of all points which a closed trapped surface passes through. A trapping boundary is defined as a connected component of the boundary of an inextendible trapped region. The concept of a marginal surface plays an important role in \cite{Ha} and also this paper, thus we give its definition here.
\begin{definition}
$\Sigma$ is called a marginal surface if on $\Sigma$, at least one of future null expansions vanishes identically, i.e.
\begin{align} \nonumber \tr \chi_{\Sigma} \equiv 0  \quad   \text{or} \quad \tr \uchi_{\Sigma} \equiv 0. \end{align}
Without loss of generality, it is assumed that $\tr \chi_{\Sigma} \equiv 0$ on $\Sigma$.
\end{definition}
With the help of double null foliations, \cite{Ha} introduces the notion of a trapping horizon as a $3$-dimensional hypersurface foliated by marginal surfaces on which $\tr \uchi \neq 0$ and $\uL \tr \chi \neq 0$, where $\tr \chi$, $\tr \uchi$, $\uL$ correspond to the ones defined in the double null foliation adapted to the marginal surfaces.

\cite{ABF1} initiates the study of isolated horizons, which is a kind of null trapping horizons with various additional conditions. Later \cite{AK1} introduces the concept of a dynamical horizon. The idea is that isolated horizons should be boundaries of black holes which are in equilibrium with their surroundings, while dynamical horizons would be for non-equilibrium states. There are a large number of papers on isolated horizons and dynamical horizons in the literature, we list some of them here: \cite{ABD}, \cite{ABF2}, \cite{AG}, \cite{AK2}, \cite{BBGB}, \cite{BF}, \cite{F}.

\cite{AG} considers the notions of a marginally trapped surface (MTS) and a marginally trapped tube (MTT). A MTS is a closed marginal surface with $\tr \uchi<0$. A MTT is a smooth $3$-dimensional hypersurface which admits a foliation by marginally trapped surfaces. The notion of a MTT relaxes the notion of a trapping horizon by dropping the condition on $\uL \tr\chi$. A spacelike MTT is a dynamical horizon. It is natural to ask whether a dynamical horizon admits different foliations by MTSs. \cite{AG} shows that the foliation by MTSs on a dynamical horizon is unique.

\cite{BF} proves an analogy of the rigidity of the foliation for a dynamical future outer trapping horizon (FOTH). According to \cite{Ha} \cite{BF}, a trapped horizon is future if $\tr \uchi<0$ while it is outer if $\uL \tr \chi<0$. \cite{BF} studies the infinitesimal deformation of a $2$-dimensional surface. Applying the infinitesimal deformation, \cite{BF} shows that, assuming that the null energy condition holds, the foliation by future outer marginally trapped surfaces (FOTSs) of a nonexpanding FOTH can be deformed by any tangential null vector field, while the foliation of a dynamical FOTH is rigid.

When considering the construction of dynamical horizons or FOTHs, a useful method is to slice the spacetime into Cauchy hypersurfaces or null hypersurfaces, then find a marginal surface in each hypersurface, and eventually glue all surfaces to form a horizon. Then the rigidity of the foliation by marginal surfaces is related to the uniqueness of the horizon constructed by the above method. Since we have the freedom to choose the way of slicing the spacetime, then if the horizons constructed from different ways of slicing are the same, it implies that the foliation by marginal surfaces isnot rigid.

This paper considers marginal surfaces and hypersurfaces foliated by marginal surfaces. We introduce a notion of a marginal tube, which is an analogy of a marginally trapped tube.
\begin{definition}
A marginal tube is a $3$-dimensional hypersurface foliated by marginal surfaces, where the corresponding null normal direction with the vanishing null expansion on the marginal surface is continuous.
\end{definition}
In the above definition, there is no assumption on the topology of marginal surfaces. The notion of a marginal tube is rather general. It includes the above mentioned trapping horizons, isolated horizons, dynamical horizons, MTTs, FOTHs. Clearly the event horizon of either a Schwarzschild or Kerr black hole is a marginal tube. A less obvious example of a marginal tube is a hyperplane (null, spacelike, or timelike) in the Minkowski spacetime. In fact, any smooth hypersurface in a Minkowski spacetime is locally a marginal tube, since its intersections with a family of parallel null hypersurfaces form a local foliation by marginal surfaces. Thus the above example implies that the notion of a marginal tube is too general for the study of black holes, if no additional condition is assumed. However it is not the goal of this paper to use the notion of a marginal tube to study black holes. We are interested in what conclusion on the foliation by marginal surfaces one can obtain with the least assumptions. A marginal tube has the least assumptions comparing to a MTT and the previous mentioned different types of horizons. 

It is natural to ask whether the foliation by marginal surfaces of a marginal tube is rigid or not. \cite{AG} \cite{BF} already answer the analogue question for dynamical horizons and FOTHs. The intuition is that the answer would depend on the type of the hypersurface. 

It is easy to show that if a marginal tube $\calH$ is null and its tangential null expansion is zero, i.e. $\tr \chi_{L_{\calH}}\equiv 0$, where $L_{\calH}$ is the tangential null vector of $\calH$, then the foliation by marginal surfaces is not rigid, and can be deformed along any tangential null vector field on $\calH$. We introduce the notion of a tangential null marginal tube to describe the above mentioned type of marginal tubes.
\begin{definition}
Let $\calH$ is a marginal tube. Suppose that $\calH$ is null and $L_{\calH}$ is the tangential future null vector field of $\calH$. $\calH$ is called a tangential null marginal tube if the future null expansion of $\calH$ with respect to $L_{\calH}$ vanishes.
\end{definition}

In fact, every spacelike section of a tangential null marginal tube is a marginal surface. In this paper, we shall prove the converse of the above assertion: if every spacelike section of a marginal tube is a marginal surface, then the marginal tube must be a tangential null marginal tube. This is the main theorem \ref{thm 5.1} in section \ref{section 5} of the paper. 

The main theorem implies that there are certain restrictions on the foliation by marginal surfaces for a spacelike marginal tube. However, theorem \ref{thm 5.1} doesnot exclude the possibility that a spacelike marginal tube admits different foliations by marginal surfaces. A simple example is the spacelike hypersurface in the Minkowski spacetime, on which the intersections with different families of parallel null hyperplanes give rise to different foliations by marginal surfaces.

As corollaries, the main theorem \ref{thm 5.1} is true if we replace marginal tubes by trapping horizons or MTTs. The theorem also implies that a dynamical horizon cannot be foliated arbitrarily by MOTSs. This is consistent with the uniqueness result in \cite{AG}.

The main tool used in the proof of theorem \ref{thm 5.1} is the double null coordinate system. We shall parametrise hypersurfaces and spacelike surfaces in the double null coordinate system. Moreover, we shall derive formulae for the geometric quantities on spacelike surfaces in terms of the parametrisation functions. These materials constitute sections \ref{section 2} - \ref{section 4}.

\section{Double null coordinate systems}\label{section 2}
\noindent
Let $(M,g)$ be a 4-dimensional time oriented spacetime. Let $C_0$ and $\uC_{0}$ be the two null hypersurfaces intersecting at a spacelike surface $\Sigma_{0,0}$. Assume that $s$ is a function on $\uC_{0}$ which decreases in the future direction, $\us$ is a function on $C_0$ which increases in future direction, and $s=\us=0$ on $\Sigma_{0,0}$. We extend $s$ and $\us$ to be optical functions in a neighbourhood of $C_0$ and $\uC_0$, i.e. $s$ and $\us$ satisfy the following eikonal equations:
\begin{align}
\nonumber
g(\nabla s, \nabla s)=g(\nabla \us, \nabla \us) =0.
\end{align}
Let $C_s$ and $\uC_{\us}$ be the corresponding level sets of $s$ and $\us$, which are all null hypersurfaces.  Denote the intersection of $C_s$ and $\uC_{\us}$ by $\Sigma_{s,\us}$. We call $\{ C_s\} \cup \{ \uC_{\us} \}$ a double null foliation in $M$.

We define the null vector fields $L', \uL'$ by
\begin{align}
\nonumber
L'=2\nabla s, \ \ \uL'=2 \nabla \us,
\end{align}
and the function $\Omega$ by
\begin{align}
\nonumber
2\Omega^{-2} = g( L',\uL').
\end{align}
Then further define
\begin{align}
\nonumber
L= \Omega^2 L', \ \ \uL=\Omega^2 \uL',
\end{align}
Then we have
\begin{align}
\nonumber
&
g(L,\uL)= 2 \Omega^2, \ \ g(L,\uL')= g( L', \uL) =2, 
\end{align}
and
\begin{align}
\nonumber
&
 L\us  = 1, \ \ \uL s =1.
\end{align}

$\{ L,\uL'\}$ and $\{L', \uL\}$ are both conjugate null frames on $\Sigma_{s,\us}$. In this paper, a conjugate null frame on a spacelike surface $\Sigma$ consists of two null normal vector field $\{L^{\Sigma}, \uL^{\Sigma} \}$ satisfying the conjugate condition $g\left(L^{\Sigma}, \uL^{\Sigma} \right) =2$.

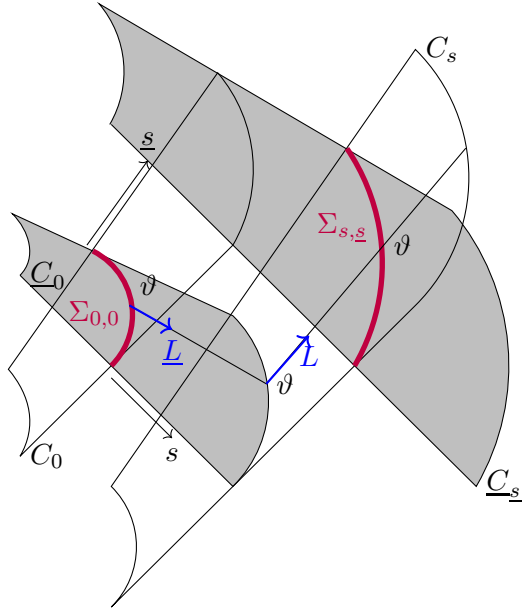
\begin{figure}[h]
\begin{center}
\begin{tikzpicture}[scale=0.8]
\draw[fill=gray!50] (-1.5,1.5) node[right] {$\uC_0$} -- (0,0) -- (2,-2) arc (-45:45:2) -- + (155:4) to [out=-45,in=45] (-1.5,1.5); 
\draw[fill=gray!50] (0,4) -- (6,-2) node[right] {$\uC_{\us}$} arc (-30:40:4) -- (-0.2,6)arc(40:-30:1.75);
\draw (-1.5,-1.5) node[right] {$C_0$} -- (0,0) -- (2,2) arc (-30:40:2.5) -- + (-125:6) to [out=-50,in=60] (-1.5,-1.5);
\draw[line width=2pt,purple] (0,0) arc (-45:0:1.19) node[left] {$\Sigma_{0,0}$} arc (0:63:1.19);
\draw[->] (-0.4,2) -- (0.6,3.4) node[above] {$\us$};
\draw[->] (0,-0.2) -- (1,-1.2) node[below] {$s$};
\draw (0,-4) -- (5,1) arc (-45:45:3) node[right] {$C_s$} -- (0,-2) arc(45:-45:1.414);
\draw[line width=2pt, purple] (4,0) arc(-30:10:3.4) node[left] {$\Sigma_{s,\us}$}arc(10:34:3.4);
\draw (0.3,1) node[above right] {$\vartheta$} -- (2.55,-0.3) node[right] {$\vartheta$} -- + (50:3) node[right] {$\vartheta$} --+ (50:5.1);
\draw[->,thick,blue] (0.3,1) -- (1,0.6) node[below] {$\uL$};
\draw[->,thick,blue] (2.55,-0.3) -- (3.25,0.5) node[below] {$L$};
\end{tikzpicture}
\end{center}
\end{figure}

We describe how to extend a coordinate system from $\Sigma_{0,0}$ to a neighbourhood of $\Sigma_{0,0}$ via the integral curves of $L$ and $\uL$:
\begin{enumerate}
\item Let $\vartheta=(\theta^1,\theta^2): U \subset \Sigma_{0,0} \rightarrow \mathbb{R}^2$ be a coordinate system on $\Sigma_{0,0}$.
\item We extend $\vartheta$ from $\Sigma_{0,0}$ to $\underline{C}_0$ such that $\vartheta$ is constant along the integral curve of $\underline{L}$.
\item We extend $\vartheta$ from $\underline{C}_0$ to a neighbourhood of $\underline{C}_0$ such that $\vartheta$ is constant along the integral curve of $L$.
\item $\{s,\underline{s},\vartheta\}$ is a coordinate system on $M$, which is called a double null coordinate system.
\end{enumerate}

By the above construction, we have
\begin{align}
\nonumber
\partial_{\us} = L.
\end{align}
We also have $\partial_s = \uL$ on $\uC_0$, but the identity may not hold outside $\uC_0$. In general $\partial_s -\uL$ is tangential to $\Sigma_{s,\us}$, therefore we define
\begin{align}
\nonumber
\vec{b} = \uL - \partial_{s}, \ \ \uL = \partial_s + \vec{b}.
\end{align}

$\{L,\uL,\partial_1= \partial_{\theta^1},\partial_2=\partial_{\theta^2}\}$ is a frame of the tangent space of $M$, and $\{ \ed L, \ed \uL, \ed \theta^1_{L,\uL}, \ed \theta^2_{L,\uL} \}$ is its dual frame. Then the metric $g$ can be expressed as
\begin{align}
\nonumber
g = 2\Omega^2 \left( \ed L \otimes \ed \uL + \ed \uL \otimes \ed L \right) + \slashg_{ab} \ed \theta^a_{L,\uL} \otimes \ed \theta^b_{L,\uL}.
\end{align}
$\{ \partial_s, \partial_{\us}, \partial_1, \partial_2 \}$ is the coordinate frame of $\{s,\us, \theta_1, \theta_2\}$ and $\{ \ed s, \ed \us, \ed \theta^1, \ed \theta^2 \}$ is the dual frame. Assume that $\vec{b}= b^1 \partial_1 + b^2 \partial_2$, then we have
\begin{align}
\nonumber
\ed L = \ed \us, \ \ 
\ed \uL = \ed s, \ \
\ed \theta^a_{L,\uL} = \ed \theta^a - b^a \ed s,
\end{align}
thus the metric in the double null coordinate system $\{s,\us, \theta_1, \theta_2\}$ takes the following form
\begin{align}
\nonumber
g = 2\Omega^2 \left( \ed s \otimes \ed \us + \ed \us \otimes \ed s \right) + \slashg_{ab} \left( \ed \theta^a - b^a \ed s \right) \otimes \left( \ed \theta^b - b^b \ed s\right).
\end{align}
We denote the Levi-Civita connection of the surface $(\Sigma_{s,\us}, \slashg)$ by $\slashnabla$ and the corresponding Christoffel symbols by $\slashGamma_{ij}^k$.

\section{Structure coefficients}\label{section 3}
\noindent
We review the structure coefficients associated with the double null foliation $\{C_s \}\cup \{\underline{C}_{\underline{s}}\}$.
\begin{definition}[Structure coeffcients]\label{SC}
Let $X,Y$ be tangential vector fields on $\Sigma_{s,\underline{s}}$. We define that
\begin{align}
&
\underline{\omega}=\underline{L} \log \Omega, 
&&
\omega= L \log \Omega, 
\nonumber
\\
\nonumber
&
\underline{\chi}(X,Y) = g(\nabla_X \underline{L},Y),  
&& 
\chi(X,Y)= g(\nabla_X L, Y), 
\\
\nonumber
&
\underline{\chi}'(X,Y) =g(\nabla_X \underline{L}',Y),  
&&
\chi'(X,Y)= g(\nabla_X L', Y),  
\\
\nonumber
&
\eta(X) = \frac{1}{2} g ( \nabla_X \underline{L}, L'),  
&&
\underline{\eta}(X) = \frac{1}{2}g(\nabla_X L, \underline{L}').
\end{align}
$\eta,\ueta$ are the torsions of the null frame $\{\uL,L'\}$ and $\{\uL' ,L\}$ respectively. $\uchi,\chi,\uchi',\chi'$ are the null second fundamental forms in the directions of $\uL,L,\uL',L'$ respectively. We can decompose them into trace and trace-free parts with respect to the metric $\slashg$,
\begin{eqnarray}
&\uchi=\hatuchi + \frac{1}{2} \tr \uchi \slashg, & \chi= \hatchi + \frac{1}{2} \tr \chi \slashg,
\nonumber
\\
\nonumber
&\uchi'=\hatuchi' + \frac{1}{2} \tr \uchi' \slashg, & \chi'= \hatchi' + \frac{1}{2} \tr \chi' \slashg.
\end{eqnarray}
The trace-free parts are called shears, and the traces are called null expansions.
\end{definition}
\begin{remark}
For any spacelike surface with a conjugate null frame, one can define the corresponding null second fundamental forms and torsions by the above definition.
\end{remark}

The null expansion $\tr \chi$ satisfies the Raychaudhuri equation
\begin{align}
\label{eqn Ray}
L \tr \chi = 2\omega \tr \chi - \frac{1}{2} |\tr \chi|^2 - |\hatchi |^2 - \R_{LL},
\end{align}
where $\R_{LL}$ is the component $\Ric(L,L)$ of the Ricci curvature.

We have the following proposition relating the covariant derivatives with the structure coefficients.
\begin{proposition}[\cite{C} Chapter 1. (1.175)]\label{SCACS} 
Let $\{e_1, e_2\}$ be a unit orthogonal frame on $\Sigma_{s,\us}$. We use the latin letters $A,B$ to denote indices varying in $\{1,2\}$.
\begin{align}
\nonumber
&\nabla_{e_A}e_B= -\frac{1}{2}\Omega^{-2}{\chi}_{AB} \underline{L} - \frac{1}{2} \Omega^{-2} \underline{\chi}_{AB} L + \slashnabla_{e_A}e_B,
\\
\nonumber
&\nabla_{e_A}\underline{L} = \eta_A \underline{L} + \underline{\chi}_{AB}  e_B,
\\
\nonumber
&\nabla_{e_A} L=-\eta_A L +{\chi}_{AB} e_B,
\\
\nonumber
&\nabla_{\underline{L}} \underline{L} =2\underline{\omega} \underline{L}, \quad \nabla_{\underline{L}} L=  -2\Omega^{2}\eta_A e_A,
\\
\nonumber
&\nabla_{L} L=2\omega L, \quad  \nabla_{L} \underline{L} = -2 \Omega^2 \underline{\eta}_A e_A.
\end{align}
\end{proposition}

\begin{remark}
From proposition \ref{SCACS}, we can derive the following equation for $\vec{b}$,
\begin{align}
\nonumber
\mathcal{L}_{L} \vec{b} =2 \Omega^2(\eta^\sharp - {\underline{\eta}}^\sharp),
\end{align}
where $\mathcal{L}$ is the Lie derivative. The derivation of the equation is simple:
\begin{align}
\nonumber
\mathcal{L}_{L} b = \left[L,\underline{L} - \partial_s \right]= \left[ L,\underline{L}\right]  - \left[ \partial_{\underline{s}},\partial_s\right] = \nabla_{L} \underline{L} -\nabla_{\underline{L}} L= 2\Omega^2(\eta^\sharp - {\underline{\eta}}^\sharp),
\end{align}
here we use the standard musical isomorphisms $\flat$ and $\sharp$ to lower and raise the indices with respect to $\slashg$.
\end{remark}

\section{Parametrisation and geometry of spacelike surfaces}\label{section 4}
\noindent
Let $\Sigma$ be a spacelike surface in $(M,g)$. If $\Sigma$ can be parametrised as the graph of two functions $f$ and $\uf$ by the following map:
\begin{align}
\vartheta \mapsto (s,\us,\vartheta)=(f(\vartheta), \uf(\vartheta), \vartheta) \in \Sigma,
\end{align}
then we call that $\Sigma$ is parametrised by $(f,\uf)$ in the double null coordinate system $\{s,\us,\vartheta\}$. For such a parametrised surface $\Sigma$, we can construct conjugate null frames on it and calculate the corresponding structure coefficients.

The tangential coordinate frame $\left\{{\dpartial_1}, {\dpartial_2} \right\}$ on $\Sigma$ is given by
\begin{align}
\nonumber
{\dpartial_i} = \partial_i + f_i \partial_s + \underline{f}_i \partial_{\underline{s}},
\end{align}
Recall that $L=\partial_{\underline{s}}, \underline{L}= \partial_s + b^i \partial_i$, thus the above is equivalent to 
\begin{align}
\nonumber
{\dpartial}_i=\left(\delta_i^j - f_i b^j \right) \partial_j + f_i \underline{L} + \underline{f}_i L = B_i^j \partial_j + f_i \underline{L} + \underline{f}_i L,
\end{align}
where 
\begin{align}
\nonumber
B_i^j= \delta_i^j - f_i b^j .
\end{align}
The intrinsic metric $\dslashg = g|_{\Sigma}$ on $\Sigma$ is given by
\begin{align}
\nonumber
\dslashg_{ij} = \dslashg\left(\dpartial_i, \dpartial_j\right) =g\left(\dpartial_i , \dpartial_j\right) = B_i^k B_j^l \slashg_{kl} + 2\Omega^2 \left(f_i \uf_j + f_j \uf_i\right).
\end{align}

In order to find a null normal frame $\left\{\duL, \dL\right\}$ on $\Sigma$, we assume that
\begin{align}
\nonumber
\left\{
\begin{aligned}
\duL= \underline{L} + \underline{\varepsilon}  L + \underline{\varepsilon}^i \partial_i,
\\
\dL= L + \varepsilon  \underline{ L} + \varepsilon^i \partial_i.
\end{aligned}
\right.
\end{align}
Then we have the following system of equations
\begin{align}
\left\{
\begin{aligned}
&g(\duL,\duL)= 4 \Omega^2 \underline{\varepsilon} + \slashg_{ij} \underline{\varepsilon}^i \underline{\varepsilon}^j=0
\\
&g(\dL , \dL)=4 \Omega^2 {\varepsilon} + \slashg_{ij}  {\varepsilon}^i  {\varepsilon}^j=0
\\
&g(\duL , {\dpartial}_i) =2\Omega^2 \underline{f}_i + 2\Omega^2 \underline{\varepsilon} f_i +  \slashg_{jk} B_i^j \underline{\varepsilon}^k=0
\\
&g(\dL , {\dpartial}_i)=2\Omega^2 {f}_i + 2\Omega^2 {\varepsilon} \underline{f}_i +  \slashg_{jk} B_i^j {\varepsilon}^k=0
\end{aligned}
\right.
\label{eqn 23}
\end{align}
Let 
\begin{align}
e^k =-2\Omega^2 f_i \left(B^{-1}\right)_j^i \left(\slashg^{-1}\right)^{jk}, \quad \underline{e}^k = -2\Omega^2 \underline{f}_i \left(B^{-1}\right)_j^i \left( \slashg^{-1} \right)^{jk},
\label{eqn 24}
\end{align}
then we can rewrite the last two equations of \eqref{eqn 23} as follows
\begin{align}
\underline{\varepsilon}^k = \underline{e}^k + \underline{\varepsilon} e^k, \quad
\varepsilon^k = e^k  + \varepsilon \underline{e}^k.
\end{align}
Substituting to the first two null condition equations of \eqref{eqn 23}, we obtain that
\begin{align}
\nonumber
\begin{aligned}
|e|^2 \underline{\varepsilon}^2 +(2 e\cdot \underline{e} + 4 \Omega^2) \underline{\varepsilon} + |\underline{e}|^2 =0,
\\
|\underline{e}|^2 {\varepsilon}^2 +(2 e\cdot \underline{e} + 4 \Omega^2) {\varepsilon} + |{e}|^2 =0,
\end{aligned}
\end{align}
where
\begin{align}
\nonumber
|e|^2 = \slashg_{ij}e^ie^j,\quad |\underline{e}|^2 = \slashg_{ij} \underline{e}^i \underline{e}^j, \quad e\cdot \underline{e} =\slashg_{ij} e^i \underline{e}^j.
\end{align}
Thus we solve the equations and obtain that
\begin{align}
\begin{aligned}
\underline{\varepsilon} = \frac{ -|\underline{e}|^2}{(2\Omega^2 + e\cdot \underline{e}) + \sqrt{(2\Omega^2 + e\cdot \underline{e})^2 -|e|^2 |\underline{e}|^2}},
\\
{\varepsilon} = \frac{ -|{e}|^2}{(2\Omega^2 + e\cdot \underline{e}) + \sqrt{(2\Omega^2 + e\cdot \underline{e})^2 -|e|^2 |\underline{e}|^2}}.
\end{aligned}
\label{eqn 28}
\end{align}
By the above calculations, we find an explicit null normal frame $\left\{ \duL, \dL \right\}$ on $\Sigma$. Associated with $\left\{ \duL, \dL \right\}$, we define the function $\dot{\Omega}$ by
\begin{align}
\nonumber
2\dot{\Omega}^2 =g\left(\duL,\dL \right)=2 \Omega^2(1+\varepsilon \cdot \underline{\varepsilon}) + \vec{\varepsilon}\cdot \vec{\underline{\varepsilon}}
\end{align}
where
\begin{align}
\nonumber
\vec{\varepsilon}\cdot \vec{\underline{\varepsilon}} =\slashg_{ij} \varepsilon^i \underline{\varepsilon}^j.
\end{align}
Then define other two null vectors $\duL',\dL'$ by
\begin{align}
\nonumber
\duL'=\dot{\Omega}^{-2} \duL, \quad \dL' = \dot{\Omega}^{-2} \dL,
\end{align}
therefore the null frames $\left\{ \duL, \dL' \right\}$ and $\left\{  \duL',\dL \right\}$ are conjugate null frames.

The structure coefficients associated with the null frame $\{ \dL, \duL\}$ are given by
\begin{align} 
&
\duchi_{ij} = g \left( \nabla_{\dpartial_i} \duL, \dpartial_j \right),
&&
\dchi_{ij} = g \left( \nabla_{\dpartial_i} \dL, \dpartial_j \right),
\nonumber
\\
&
\deta_i = \frac{1}{2} g \left( \nabla_{\dpartial_i} \duL, \dL' \right),
&&
\dot{\ueta}_i = \frac{1}{2} g \left( \nabla_{\dpartial_i} \dL, \duL' \right).
\nonumber
\end{align}
The following proposition gives the formulae of $\duchi,\dchi,\deta$.
\begin{proposition}\label{pro 4.1}
Let $\Sigma$ be a spacelike surface parametrised by $\{ f, \uf\}$ in the double null coordinate system $\{s,\us,\vartheta\}$. Then the structure coefficients of $\Sigma$ have the following formulae.
\begin{align}
\duchi_{ij}
=&
\uchi_{ij} + \uvarepsilon \chi_{ij} + (\vec{\uvarepsilon}\cdot b -2\Omega^2 \uvarepsilon) \slashnabla^2_{ij} f -2\Omega^2 \slashnabla_{ij}^2\underline{f}
-2 \sym \left\{ 
\left[ \chi(\vec{\uvarepsilon}) + 2\Omega^2 \ueta \right] \otimes \dslashd \uf
\right\}_{ij}
\nonumber
\\
\nonumber
&
+2\sym \left\{ 
\left[ \slashnabla b \cdot \vec{\uvarepsilon} - \uchi(\uvarepsilon) - \uvarepsilon \uchi(b) -2\Omega^2 \uvarepsilon \eta \right] \otimes \dslashd f 
\right\}_{ij}
\\
\nonumber
&
+ 2 \left[  2\Omega^2 \eta\cdot \vec{\uvarepsilon} +\chi(b,\vec{\uvarepsilon}) + 2\Omega^2 b\cdot \eta  \right] \sym(\dslashd f \otimes \dslashd \uf)_{ij} 
-4\omega \Omega^2 (\dslashd \uf \otimes \dslashd \uf)_{ij}
\\
\nonumber
&
+ \left[
2\uchi(b,\vec{\uvarepsilon}) + \uchi(b,b) + \uvarepsilon \chi(b,b) +4 \Omega^2 \uvarepsilon b\cdot \eta - \slashnabla_b b\cdot \vec{\uvarepsilon} +\frac{\partial b}{\partial s}\cdot \vec{\uvarepsilon} + 4 \Omega^2 \uvarepsilon \uomega
\right] 
\\
&\ \ \ \ \cdot
\left( \dslashd f \otimes \dslashd f \right)_{ij},
\end{align}

\begin{align} 
\dchi_{ij}
=&
\chi_{ij} + \varepsilon \uchi_{ij} + (\vec{\varepsilon}\cdot b -2\Omega^2 ) \slashnabla^2_{ij} f -2\Omega^2 \varepsilon \slashnabla_{ij}\underline{f}
\nonumber
\\
\nonumber
&
-2 \sym \left\{
\left[\chi(\vec{\varepsilon}) + 2\Omega^2 \varepsilon \ueta \right] \otimes \dslashd \uf
\right\}_{ij}
\\
&
+2\sym \left\{
\left[ \slashnabla b \cdot \vec{\varepsilon} - \uchi(\vec{\varepsilon}) - \varepsilon \uchi(b) -\chi(b) -2\Omega^2  \eta \right] \otimes \dslashd f
\right\}_{ij}
\nonumber
\\
&
+2 \left[  2\Omega^2 \eta\cdot \vec{\varepsilon} +\chi(b,\vec{\varepsilon}) + 2\Omega^2 \varepsilon b\cdot \eta  \right] \sym(\dslashd f \otimes \dslashd \uf)_{ij} 
-4\varepsilon \omega \Omega^2 (\dslashd \uf \otimes \dslashd \uf)_{ij}
\nonumber
\\
&
+ \left[
2\uchi(b,\vec{\varepsilon}) + \chi(b,b) + \varepsilon \uchi(b,b) +4 \Omega^2 b\cdot \eta - \slashnabla_b b\cdot \vec{\varepsilon} +\frac{\partial b}{\partial s}\cdot \vec{\varepsilon} + 4 \Omega^2  \uomega
\right] 
\nonumber
\\
&\ \ \ \ \cdot
\left( \dslashd f \otimes \dslashd f \right)_{ij},
\end{align}

\begin{align} 
\deta_i
=&
\left[ 2\Omega^2 (1 +\varepsilon \uvarepsilon) + \vec{\varepsilon} \cdot \vec{\uvarepsilon} \right]^{-1}
\cdot
\nonumber
\\
&
\left\{
2\Omega^2 \eta_i + 2 \Omega^2 \varepsilon \uvarepsilon \ueta_i - \chi(\vec{\uvarepsilon})_i + \uvarepsilon \chi( \vec{\varepsilon})_i +\uchi( \vec{\varepsilon})_i - \varepsilon \uchi(\vec{\uvarepsilon})_i
\right.
\nonumber
\\
& \quad
+
f_i
\left[
4\Omega^2 \uomega -2\Omega^2 b\cdot \eta -2 \Omega^2 \varepsilon \uvarepsilon b\cdot \eta + 2\Omega^2 \vec{\uvarepsilon}\cdot \eta -2\Omega^2 \uvarepsilon \eta\cdot \vec{\varepsilon}
\right.
\nonumber
\\
&\qquad \qquad
\left.
+\chi(b,\vec{\uvarepsilon}) - \uvarepsilon \chi (b, \vec{\varepsilon}) + \uchi(\vec{\uvarepsilon},\vec{\varepsilon}) - \uchi(b,\vec{\varepsilon}) + \varepsilon \uchi(b,\vec{\uvarepsilon}) -(\slashnabla_{\vec{\uvarepsilon}} b) \cdot \vec{\varepsilon}
\right]
\nonumber
\\
&\quad
\left. 
+ \uf_i \left[
-2\Omega^2 \ueta\cdot \vec{\varepsilon} + 2\Omega^2 \ueta \cdot \vec{\uvarepsilon} + 4\Omega^2 \varepsilon \uvarepsilon \omega + \chi(\vec{\uvarepsilon}, \vec{\uvarepsilon})
\right]
+
2\Omega^2 \varepsilon (\dslashd \uvarepsilon)_i + (\slashnabla_i \vec{\uvarepsilon})\cdot \vec{\varepsilon}
\right\}.
\end{align}
In the above formulae, $\dslashd$ in $\dslashd f, \dslashd \uf, \dslashd \uvareps$ is the differential operator on $\Sigma$. $\slashnabla$ in $ \slashnabla^2_{ij} f$, $\slashnabla^2_{ij} \uf$, $\slashnabla_i \vec{\uvareps} $
is the pull back of the covariant derivative $\slashnabla$ via the embedding of $\Sigma$, i.e. for any vector field $V$ on $\Sigma$,
\begin{align}
\nonumber
\slashnabla_i V^k = \dpartial_i V^k + \slashGamma_{ij}^k V^j, 
\end{align}
and more generally for any tensor field $T$ on $\Sigma$,
\begin{align}
\nonumber
\slashnabla_{i} T_{i_1 \cdots i_k}^{j_1 \cdots j_l}
=
\dpartial_i T_{i_1 \cdots i_k}^{j_1 \cdots j_l} 
-  \slashGamma_{i i_m}^{r}  T_{i_1 \cdots \underset{\hat{i_m}}{r}\cdots i_k}^{j_1 \cdots j_l}
+ \slashGamma_{i s}^{j_n} T_{i_1 \cdots i_k}^{j_1 \cdots \overset{\hat{j_n}}{s}\cdots  j_l}.
\end{align}
The inner products $\cdot$ in $\vec{\uvareps}\cdot b$, $\slashnabla b \cdot \vec{\uvarepsilon}$, $\slashnabla_b b\cdot \vec{\uvarepsilon}$, $\frac{\partial b}{\partial s}\cdot \vec{\uvarepsilon}$,  $\vec{\varepsilon}\cdot b$, $\slashnabla b \cdot \vec{\varepsilon}$, $\slashnabla_b b\cdot \vec{\varepsilon}$, $\vec{\varepsilon} \cdot \vec{\uvarepsilon}$, $(\slashnabla_i \vec{\uvarepsilon})\cdot \vec{\varepsilon}$ are taken with respect to the metric $\slashg$.
\end{proposition}

We give the proof of the above proposition in the appendix. For the rest of this section, we consider two special cases of proposition \ref{pro 4.1}.

The first case is that $\Sigma$ is embedded in the null hypersurface $C_{s=s_0}$. This case was also investigated in \cite{KLR} \cite{L} \cite{An} before.

In this case, the parametrisation function $f$ of $\Sigma$ is constant, $f\equiv s_0$. $\Sigma$ is parametrised by the following map:
\begin{align}
\nonumber
\vartheta \mapsto (s,\us,\vartheta) = (s_0, \uf(\vartheta), \vartheta).
\end{align}
The tangential coordinate frame $\{\dpartial_1, \dpartial_2\}$ of $\Sigma$ is
\begin{align}
\nonumber
\dpartial_i=\partial_i + \uf_i \partial_{\us} = \partial_i + \uf_i L.
\end{align}
The intrinsic metric $\dslashg= g|_{\Sigma}$ on $\Sigma$ is
\begin{align}
\nonumber
\dslashg_{ij}=\slashg_{ij}.
\end{align}
The null normal frame $\left\{ \duL,\dL \right\}$ is
\begin{align}
\nonumber
\left\{
\begin{aligned}
& \duL= \uL + \uvareps L + \uvareps^i \partial_i,
\\
& \dL=L,
\end{aligned}
\right.
\end{align}
where
\begin{align}
\nonumber
\uvareps = - \Omega^2 \left( \slashg^{-1} \right)^{ij} \uf_i \uf_j = - \Omega^2 \leve \dslashd \uf \rive_{\slashg}^2,  \quad \uvareps^k= -2\Omega^2 \uf_i \left( \slashg^{-1} \right)^{ik}.
\end{align}
The inner product of $\duL$ and $\dL$ is
\begin{align}
\nonumber
2\dot{\Omega}^2= g \left( \duL,\dL \right)= 2\Omega^2 \Rightarrow \dot{\Omega}=\Omega,
\end{align}
thus
\begin{align}
\nonumber
\dL'=\dot{\Omega}^{-2} \dL = \Omega^{-2} L = L',
\end{align}
and $\left\{ \duL,\dL'=L' \right\}$ is a conjugate null frame on $\Sigma$.

By proposition \ref{pro 4.1}, we obtain the structure coefficients of $\Sigma$ in the first case as follows.
\begin{proposition}\label{pro 4.2}
Let $\Sigma$ be a spacelike surface embedded in the null hypersurface $C_{s=s_0}$. Assume that $\Sigma$ is parametrised by $(s_0, \uf)$ in the double null coordinate system $(s, \us, \vartheta)$. Then the structure coefficients of $\Sigma$ have the following formulae.
\begin{align}
&
\begin{aligned}
\duchi_{ij}
=&
\uchi_{ij} -\Omega^2 \leve \dslashd \uf \rive^2_{\slashg} \chi_{ij}  -2\Omega^2 \slashnabla_{ij}^2\underline{f} -4\Omega^2 \sym \left\{\ueta \otimes \dslashd \uf \right\}_{ij} 
- 4\omega\Omega^2 (\dslashd \uf \otimes \dslashd \uf )_{ij} 
\nonumber
\\
&
+ 4\Omega^2 \sym \left\{\dslashd\uf \otimes \chi(\slashnabla \uf) \right\}_{ij},
\nonumber
\\
\dtr \duchi 
= & \tr \uchi - 2\Omega^2 \slashDelta \uf - \Omega^2 \left|\dslashd \uf \right|_{\slashg}^2 \tr \chi -4\Omega^2 \ueta \cdot \dslashd \uf - 4 \Omega^2 \omega \left| \dslashd \uf \right|^2_{\slashg} + 4 \Omega^2 \chi ( \slashnabla \uf, \slashnabla \uf),
\nonumber
\end{aligned}
\\
&
\begin{aligned}
\dchi_{ij} = \chi_{ij},  \quad \dtr \dchi = \tr \chi, \quad \deta_i = \eta_i + \chi (\slashnabla \uf )_i.
\nonumber
\end{aligned}
\end{align}
In the above formulae, $\slashnabla^i \uf = \left( \slashg^{-1} \right)^{ij} \uf_j$, $\slashDelta \uf = \left( \slashg^{-1} \right)^{ij} \slashnabla_{ij}^2 \uf$. The inner product $\cdot$ in $\ueta \cdot \dslashd \uf $ is taken with respect to $\slashg$. $\dtr$ in $\dtr \dchi$ is the trace operator with respect to $\dslashg$, although numerically the same as $\slashg$.
\end{proposition}

The second case is that $\Sigma$ is embedded in the null hypersurface $\uC_{\us=\us_0}$. In this case, the parametrisation function $\uf$ of $\Sigma$ is constant, $\uf \equiv \us_0$. $\Sigma$ is parametrised by the following map:
\begin{align}
\nonumber
\vartheta \mapsto (s,\us, \vartheta) = (f(\vartheta), \us_0, \vartheta).
\end{align}
The tangential coordinate frame $\left\{ \dpartial_1,\dpartial_2\right\}$ of $\Sigma$ is
\begin{align}
\nonumber
\dpartial_i = \partial_i + f_i \partial_s = \partial_i + f_i (\uL-b^j \partial_j) = B_i^j \partial_j +f_i \uL,
\end{align}
where
\begin{align}
\nonumber
B_i^j=\delta_i^j -f_i b^j.
\end{align}
The intrinsic metric $\dslashg= g|_{\Sigma}$ is
\begin{align}
\nonumber
\dslashg_{ij}=B_i^k B_j^l \slashg_{kl} = \slashg_{ij} - ( \slashg_{il} b^l f_j + \slashg_{jl} b^l f_i) + f_i f_j |b|^2_{\slashg}.
\end{align}
The null normal frame $\left\{ \duL,\dL \right\}$ on $\Sigma$ is
\begin{align}
\nonumber
\left\{
\begin{aligned}
\duL&=\uL,
\\
\dL&=L+ \varepsilon \uL +  \varepsilon^i \partial_i,
\end{aligned}
\right.
\end{align}
where
\begin{align}
\nonumber
\vareps=-\Omega^2 \dslashg^{ij}f_i f_j = - \Omega^2 \slashg^{kl} {(B^{-1})}_k^i {(B^{-1})}_l^j f_i f_j, \quad  \vareps^k=- 2\Omega^2 \slashg^{kl} {(B^{-1})}_l^i f_i.
\end{align}
The inner product of $\duL$ and $\dL$ is
\begin{align}
\nonumber
2\dot{\Omega}^2= g \left( \duL,\dL \right)= 2\Omega^2 \Rightarrow \dot{\Omega}=\Omega,
\end{align}
thus
\begin{align}
\nonumber
\dL' = \dot{\Omega}^{-2} \dL = \Omega^{-2} \left( L+ \varepsilon \uL +  \varepsilon^i \partial_i \right) = L'+ \varepsilon' \uL +  \varepsilon'^i \partial_i,
\end{align}
where
\begin{align}
\nonumber
\vareps' = \Omega^{-2} \vareps= - \slashg^{kl} {(B^{-1})}_k^i {(B^{-1})}_l^j f_i f_j, \quad  \vareps'^{k}= \Omega^{-2} \vareps^k=- 2 \slashg^{kl} {(B^{-1})}_l^i f_i.
\end{align}
$\left\{\duL, \dL' \right\}$ is a conjugate null frame on $\Sigma$. 

By proposition \ref{pro 4.1}, we obtain the structure coefficients on $\Sigma$ associated with the conjugate null frame $\left\{ \duL,\dL' \right\}$ in the second case as follows.
\begin{proposition}\label{pro 4.3}
Let $\Sigma$ be a spacelike surface embedded in the null hypersurface $\uC_{\us=\us_0}$. Assume that $\Sigma$ is parametrised by $(f,\us_0)$ in the double null coordinate system $(s, \us, \vartheta)$. Then the structure coefficients of $\Sigma$ have the following formulae.
\begin{align}
\duchi_{ij}=&\uchi_{ij} - 2 \sym \left\{  \uchi(b)  \otimes \dslashd f \right\}_{ij} + \uchi(b,b) f_i f_j,
\nonumber
\\
\dchi'_{ij} 
=&
\chi'_{ij} + \vareps' \uchi_{ij} + ( b\cdot \vec{\vareps'} -2) \slashnabla^2_{ij} f
\nonumber
\\
&
+2\sym \left\{  \left[ \slashnabla b\cdot \vec{\vareps'} - \uchi(\vec{\vareps'}) -\vareps' \uchi(b) - \chi'(b) -2 \eta \right] \otimes \dslashd f
\right\}_{ij}
\nonumber
\\
&
+ \left[
2\uchi(b,\vec{\vareps'}) + \vareps' \uchi(b,b) + \chi'(b,b) + 4 b\cdot \eta - \slashnabla_b b \cdot \vec{\vareps'}
-\partial_s b\cdot \vec{\vareps'} - 4  \uomega
\right] f_i f_j,
\nonumber
\\
\deta_i =& \eta_i + \frac{1}{2} \uchi(\vec{\vareps'})_i + \left[ 2\uomega-b\cdot \eta -\frac{1}{2} \uchi(b,\vec{\vareps'}) \right] f_i.
\nonumber
\end{align}
In the above formulae, the inner products $\cdot$ in $b\cdot \vec{\vareps'}$, $\slashnabla b \cdot \vec{\vareps'}$, $b\cdot \eta$, $\slashnabla_b b \cdot \vec{\vareps'}$, $\partial_s b \cdot \vec{\vareps'}$ are taken with respect to the metric $\slashg$. 
\end{proposition}

\section{The main result}\label{section 5}
\noindent
We state and prove the main result of the paper.

\begin{theorem}\label{thm 5.1}
Let $(M,g)$ be a 4-dimensional spacetime and $\{s,\us,\vartheta\}$ be a double null coordinate system on $(M,g)$. Let $\calH$ be a smooth hypersurface containing $\Sigma_{0,0}$. Assume that $\calH$ can be parametrised as the graph of a function $h$ of independent variables $(\us,\vartheta)$,
\begin{align}
\nonumber
&
\calH= \left\{ (s,\us,\vartheta) = \left( h(\us,\vartheta), \us, \vartheta \right) \right\},
\\
\nonumber
&
\psi_h: (\us,\vartheta) \mapsto  (s,\us,\vartheta) = \left( h(\us,\vartheta), \us, \vartheta \right) \in \calH.
\end{align}
If every spacelike section $\Sigma$ of $\calH$ is a marginal surface with 
\begin{align}
\nonumber
\tr \chi_{\Sigma} \equiv 0,
\end{align}
then $\calH$ is null which is equivalent to $h\equiv 0$, therefore $\calH$ is a tangential null marginal tube.
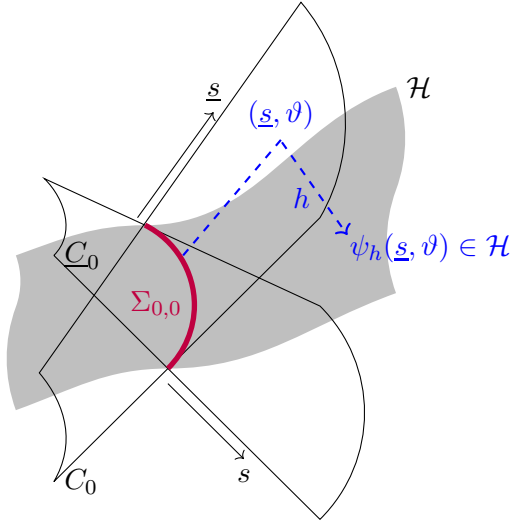
\begin{center}
\begin{tikzpicture}[scale=1]
\path[fill=gray!50] 
(-2,1.5) arc(20:-20:1) arc(160:200:2) 
[out=20,in=180] to (0,0)
[out=0,in=-160] to (3,1) 
arc(200:160:2) arc(-20:20:2) node[right] {$\calH$}
[out=-160,in=0] to (-0.3,1.9)
[out=180,in=20] to (-2,1.5);
\draw (-1.5,1.5) node[right] {$\uC_0$} -- (0,0) -- (2,-2) arc (-45:45:2) -- + (155:4) to [out=-45,in=45] (-1.5,1.5); 
\draw (-1.5,-1.5) node[right] {$C_0$} -- (0,0) -- (2,2) arc (-30:40:2.5) -- + (-125:6) to [out=-50,in=60] (-1.5,-1.5);
\draw[line width=2pt,purple] (0,0) arc (-45:0:1.19) node[left] {$\Sigma_{0,0}$} arc (0:63:1.19);
\draw[->] (-0.4,2) -- (0.6,3.4) node[above] {$\us$};
\draw[->] (0,-0.2) -- (1,-1.2) node[below] {$s$};
\draw[thick,blue,dashed] (0.2,1.5) --+(50:2) node[above] {$(\us,\vartheta)$};
\draw[thick,blue,dashed,->] (1.5,3) --+ (-55:0.9) node[left] {$h$} --+ (-55:1.3) node[below right] {$\psi_{h}(\us, \vartheta) \in \calH$} --+ (-55:1.5);
\end{tikzpicture}
\end{center}
\end{theorem}

\begin{proof}[Proof of theorem \ref{thm 5.1}]
The goal is to show that $h\equiv 0$. Consider the set $N$ of points $(\us,\vartheta)$ where $\partial_{\us} h \neq 0$. $N$ is an open set. We shall show that $N$ is empty by contradiction. 

If $N$ is not empty, choose a point $(\us_0, \vartheta_0)\in N$. By the implicit function theorem, there is a neighbourhood $U_0$ of the point $(h(\us_0,\vartheta_0) ,\us_0,\vartheta_0 )$ in $\calH$, such that $\calH$ can be parametrised by $(s,\vartheta)$ in $U_0$.

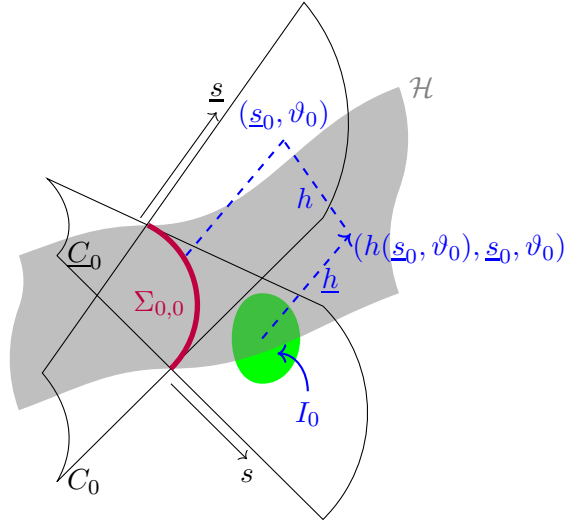
\begin{center}
\begin{tikzpicture}[scale=1]
\path[fill=green] (1.2,-0.2)
to [out=0,in=-90] (1.7,0.4)
to [out=90,in=0] (1.2,1)
to [out=180,in=90] (0.8,0.4)
to [out=-90,in=180] cycle; 
\path[fill=gray,opacity=0.5] 
(-2,1.5) arc(20:-20:1) arc(160:200:2) 
[out=20,in=180] to (0,0)
[out=0,in=-160] to (3,1) 
arc(200:160:2) arc(-20:20:2) node[right] {$\calH$}
[out=-160,in=0] to (-0.3,1.9)
[out=180,in=20] to (-2,1.5);
\draw (-1.5,1.5) node[right] {$\uC_0$} -- (0,0) -- (2,-2) arc (-45:45:2) -- + (155:4) to [out=-45,in=45] (-1.5,1.5); 
\draw (-1.5,-1.5) node[right] {$C_0$} -- (0,0) -- (2,2) arc (-30:40:2.5) -- + (-125:6) to [out=-50,in=60] (-1.5,-1.5);
\draw[line width=2pt,purple] (0,0) arc (-45:0:1.19) node[left] {$\Sigma_{0,0}$} arc (0:63:1.19);
\draw[->] (-0.4,2) -- (0.6,3.4) node[above] {$\us$};
\draw[->] (0,-0.2) -- (1,-1.2) node[below] {$s$};
\draw[thick,blue,dashed] (0.2,1.5) --+(50:2) node[above] {$(\us_0,\vartheta_0)$};
\draw[thick,blue,dashed] (1.5,3) --+ (-55:0.9) node[left] {$h$} --+ (-55:1.3) node[below right] {$(h(\us_0,\vartheta_0), \us_0, \vartheta_0)$} --+ (-55:1.5);
\draw[thick,blue,dashed,->] (1.2,0.4) --+ (50:1) node[right] {$\underline{h}$} --+ (50:1.8);
\draw[thick,blue,->] (1.8,-0.3) node[below] {$I_0$}
to [out=90,in=0] (1.4,0.2);
\end{tikzpicture}
\end{center}

More precisely, there exist a neighbourhood $U_0$ of the point $(h(\us_0,\vartheta_0) ,\us_0,\vartheta_0 )$ in $\calH$, a neighbourhood $I_0$ of the point $(h(\us_0,\vartheta_0) ,\vartheta_0)$ in the $(s,\vartheta)$ domain, and a function $\underline{h}(s,\vartheta)$ such that $U_0$ can be parametrised by the following map
\begin{align}
\nonumber
&
\psi_{\underline{h}}: 
I_0 \rightarrow U_0, 
\quad 
(s,\vartheta) \mapsto  ( s, \underline{h}( s,\vartheta), \vartheta )
\\
\nonumber
&
U_0= \left\{(s,\us,\vartheta) = ( s, \underline{h}( s,\vartheta), \vartheta ), (s,\vartheta) \in I_0\right\}.
\end{align}

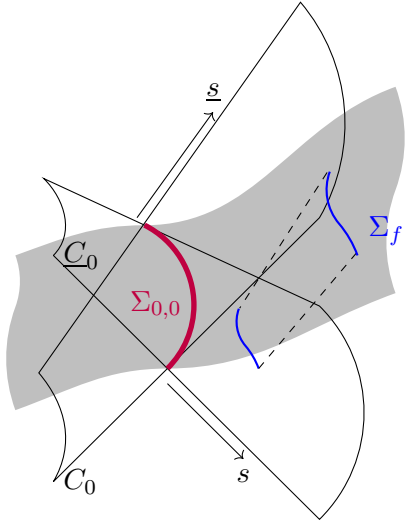
\begin{center}
\begin{tikzpicture}[scale=1]
\path[fill=gray!50] 
(-2,1.5) arc(20:-20:1) arc(160:200:2) 
[out=20,in=180] to (0,0)
[out=0,in=-160] to (3,1)
arc(200:160:2) arc(-20:20:2)
[out=-160,in=0] to (-0.3,1.9)
[out=180,in=20] to (-2,1.5);
\draw (-1.5,1.5) node[right] {$\uC_0$} -- (0,0) -- (2,-2) arc (-45:45:2) -- + (155:4) to [out=-45,in=45] (-1.5,1.5); 
\draw (-1.5,-1.5) node[right] {$C_0$} -- (0,0) -- (2,2) arc (-30:40:2.5) -- + (-125:6) to [out=-50,in=60] (-1.5,-1.5);
\draw[line width=2pt,purple] (0,0) arc (-45:0:1.19) node[left] {$\Sigma_{0,0}$} arc (0:63:1.19);
\draw[->] (-0.4,2) -- (0.6,3.4) node[above] {$\us$};
\draw[->] (0,-0.2) -- (1,-1.2) node[below] {$s$};
\draw[thick,blue] (1.2,0) arc(20:40:1) arc(220:160:0.5);
\draw[dashed] (1.2,0) --+(50:2);
\draw[thick,blue] (2.5,1.5) node[above right] {$\Sigma_f$}arc(20:40:1.4) arc(220:160:0.7);
\draw[dashed] (2.1,2.6) -- (0.95,0.8);
\end{tikzpicture}
\end{center}

Consider the spacelike surface $\Sigma_{f}$ in $U_0$ parametrised by the following map
\begin{align}
\nonumber
&
\varphi_f:\;
\vartheta \;
\mapsto \;
 (f(\vartheta), \vartheta) \;
\overset{\psi_{\underline{h}}}{\mapsto} \;
\left( f(\vartheta),\underline{h}(f(\vartheta) ,\vartheta), \vartheta \right),
\\
\nonumber
&
\Sigma_f=\left\{ (s,\us,\vartheta) = \left( f(\vartheta),\underline{h}(f(\vartheta) ,\vartheta), \vartheta \right), (s,\vartheta) \in I_0 \right\},
\end{align}
where $f$ is a function of independent variable $\vartheta$. 

Define $\uf(\vartheta) = \underline{h}(f(\vartheta), \vartheta)$, then $\Sigma_{f}$ is parametrised by $(f,\uf)$ in the double null coordinate system $\{s, \us, \vartheta \}$.
\begin{align}
\nonumber
\Sigma_f= \left\{ (s,\us,\vartheta) = \left( f(\vartheta), \uf(\vartheta) ,\vartheta  \right) \right\}.
\end{align}
We have the freedom to choose $f$. By the assumption, any $\Sigma_f$ is a marginal surface. In the following, we consider two choices of the functions $f$.
\begin{enumerate}
\item Choose $f\equiv s_0$ being constant. The corresponding surface is $\Sigma_{f\equiv s_0}$. 
\item Choose $\bar{f}$ such that $\bar{f}(\vartheta_0) = s_0$ and $\partial_{\vartheta} \bar{f} (\vartheta_0)=0$. The corresponding surface is $\Sigma_{\bar{f}}$.
\end{enumerate}

We shall use the assumption that the null expansions of both $\Sigma_{f\equiv s_0}$ and $\Sigma_{\bar{f}}$ are zero to derive a contradiction.

\begin{enumerate}
\item
By proposition \ref{pro 4.2}, the background null expansion $\tr \chi$ with respect to $L$ in the double null coordinate system $\{s, \us, \vartheta \}$ vanishes on $\Sigma_{f\equiv s_0}$,
\begin{align}
\nonumber
\tr \chi =0 \text{ on } \Sigma_{f\equiv s_0}.
\end{align}

\item
Let $p_0$ be the point $(s_0,\underline{h}_0(s_0,\vartheta_0),\vartheta_0)$ in $\Sigma_{f\equiv s_0}$. Since $\bar{f}(\vartheta_0) =s_0$ and $\partial_{\vartheta}\bar{f} (\vartheta_0) =0$, then the spacelike surface $\Sigma_{\bar{f}}$ is tangential to $\Sigma_{f\equiv s_0}$ at the point $p_0$, i.e. they have the same tangent space and normal space at $p_0$. 

$\Sigma_{\bar{f}}$ is parametrised by $(\bar{f},\bar{\uf})$ where $\bar{\uf}(\vartheta) = \underline{h} ( \bar{f}(\vartheta),\vartheta)$. Since $\partial_{\vartheta} \bar{f}(\vartheta_0)=0$, by the equations \eqref{eqn 24} - \eqref{eqn 28}, the corresponding $e^k, \vareps^k,\vareps$ of $\Sigma_{\bar{f}}$ are all equal to zero at $p_0$. The intrinsic metric $\bslashg$ of $\Sigma_{\bar{f}}$ at $p_0$ is equal to $\slashg$. Thus applying proposition \ref{pro 4.1} to $\Sigma_{\bar{f}}$ at $p_0$, we obtain that the null expansion $\btr \bchi$ of $\Sigma_{\bar{f}}$ at $p_0$ is
\begin{align}
\nonumber
\btr \bchi (p_0) =\tr \chi(p_0) - 2 \Omega^2 \left( \slashg^{-1}\right)^{ij}\slashnabla^2_{ij}  \bar{f}(\vartheta_0).
\end{align}
The assumption that $\btr \bchi(p_0)=0$ implies that $\left( \slashg^{-1}\right)^{ij}\slashnabla^2_{ij}  \bar{f}(\vartheta_0)=0$.
\end{enumerate}

Therefore we show that every $\bar{f}$, with $\bar{f}(\vartheta_0)=s_0, \partial_{\vartheta} \bar{f}(\vartheta_0) =0$, must satisfy the equation
\begin{align}
\nonumber
\left( \slashg^{-1}\right)^{ij}\slashnabla^2_{ij}  \bar{f}(\vartheta_0)=0.
\end{align}
This is impossible, thus we conclude that the set $N$ where $\partial_{\us} h\neq 0$ is empty. 

Therefore $h\equiv 0$ at least in a small neighbourhood of $\Sigma_{0,0}$, thus $\calH$ is null in a neighbourhood of $\Sigma_{0,0}$. To show that $\calH$ is null everywhere, we just apply the local result to any spacelike section of $\calH$. Note that the double null coordinate system can be constructed in a neighbourhood of any smooth spacelike surface, hence we can use the above local result to any smooth spacelike section of $\calH$. Therefore $\calH$ is null and the theorem is proved.
\end{proof}

\begin{remark}
Note that the proof only makes use of the formula for $\dchi_{ij}$ in proposition \ref{pro 4.1}. Since proposition \ref{pro 4.1} is proved in arbitrary $4$-dimensional Lorentzian manifolds, then theorem \ref{thm 5.1} is also valid for arbitrary $4$-dimensional spacetimes. It applies to spacetimes without any energy condition, or in alternative theories of gravity.
\end{remark}

We call the spacetime satisfying the null convergence condition if $\R(X,X)\geq 0$ for any null vector $X$, see \cite{HE}. For such a spacetime, we show that the tangential null marginal tube is shear-free. This is a direct corollary of the Raychaudhuri equation \eqref{eqn Ray}, since
\begin{align}
\nonumber
0 = L \tr \chi = 2\omega \tr \chi - \frac{1}{2} |\tr \chi|^2 - |\hatchi |^2 - \R_{LL} \leq - |\hatchi|^2
\quad
\Rightarrow
\quad
| \hatchi |^2 =0.
\end{align}

\appendix

\section{Proof of proposition \ref{pro 4.1}}
\begin{lemma}\label{lem A.1}
Let $\{s, \us, \vartheta\}$ be a double null coordinate system of $(M,g)$ constructed as in section \ref{section 2}. We have
\begin{align}
\nonumber
&
\nabla_{\partial_i} \partial_s = \left( \uchi_i^{\phantom{i}j} - \slashnabla_i b^j \right) \partial_j + \left( \eta_i + \frac{1}{2} \Omega^{-2} \chi(b)_i \right) \uL + \frac{1}{2} \Omega^{-2} \uchi(b)_i L,
\\
\nonumber
&
\nabla_{\partial_i} \partial_{\us} = \eta_i L + \chi_i^{\phantom{i}j} \partial_j,
\\
\nonumber
&
\nabla_{\partial_s} \partial_{\us} = -2 \Omega^2 \eta^{\sharp} - \chi(b)^{\sharp} - \ueta(b) \cdot L,
\\
\nonumber
&
\nabla_{\partial_{\us}} \partial_{\us} = 2\omega L,
\\
\nonumber
&
\nabla_{\partial_s} \partial_s = \left[ 2\uomega - 2 \eta(b) - \frac{1}{2} \Omega^{-2} \chi(b,b) \right] \uL - \frac{1}{2} \Omega^{-2} \uchi(b,b) L - 2 \uchi(b)^{\sharp} + \slashnabla_{b} b - \partial_s b.
\end{align}
\end{lemma}
\begin{proof}
It is an immediate corollary of proposition \ref{SCACS} and $\partial_s=\uL-\vec{b}$.
\end{proof}

\begin{lemma}\label{lem A.2}
Let $\Sigma$ be a spacelike surface parametrised by $(f,\uf)$ in the double null coordinate system $\{s,\us,\vartheta\}$. $\{ \dpartial_1 , \dpartial_2\}$ is the tangential coordinate vector frame on $\Sigma$. We have
\begin{align}
\nonumber
\nabla_{\dpartial_i} \dpartial_j = \Pi_{ij}^k \partial_k + \Pi_{ij}^L L + \Pi_{ij}^{\uL} \uL,
\end{align}
where
\begin{align}
\nonumber
\Pi_{ij}^k 
=&
\slashGamma_{ij}^{\phantom{ij}k} - f_{ij} b^k + 2\sym(  \dslashd f \otimes \uchi )_{ij}^{\phantom{ij}k} - 2\sym ( \dslashd f \otimes \slashnabla b  )_{ij}^{\phantom{ij}k} + 2 \sym ( \dslashd \uf \otimes \chi  )_{ij}^{\phantom{ij}k} 
\\
&
+ ( \dslashd f \otimes \dslashd f )_{ij} \cdot \left[ \slashnabla_b b - 2 \uchi(b) - \partial_s b \right]^k,
\nonumber
\\
\nonumber
\Pi_{ij}^L 
=&
\uf_{ij} - \frac{1}{2}\Omega^{-2} \uchi_{ij} + \Omega^{-2} \sym ( \dslashd f \otimes \uchi(b) )_{ij} + 2\sym( \dslashd \uf \otimes \ueta )_{ij}
\\
&
-\frac{1}{2} \Omega^{-2} \uchi(b,b) ( \dslashd f \otimes \dslashd f)_{ij} + 2 \omega ( \dslashd \uf \otimes \dslashd \uf)_{ij} - 2 \eta(b) \cdot\sym(\dslashd f \otimes \dslashd \uf)_{ij},
\nonumber
\\
\nonumber
\Pi_{ij}^{\uL}
=&
f_{ij} - \frac{1}{2} \Omega^{-2} \chi_{ij} + 2 \sym( \dslashd f \otimes \eta)_{ij} + \Omega^{-2} \sym(\dslashd f \otimes \chi(b) )_{ij}
\\
&
+ \left[ 2\uomega - 2 \eta(b) - \frac{1}{2} \Omega^{-2} \chi(b,b) \right] (\dslashd f \otimes \dslashd f)_{ij},
\nonumber
\end{align}
where
\begin{align}
\nonumber
&
\sym( P \otimes Q)_{ij} = \frac{1}{2} (P_i Q_j + P_j Q_i),
\\
\nonumber
&
\sym (P \otimes T )_{ij}^{\phantom{ij}k} = \frac{1}{2} ( P_i T_j^k + P_j T_i^k ).
\end{align}
\end{lemma}
\begin{proof}
\begin{align}
\nonumber
\nabla_{\dpartial_i} \dpartial_j 
=&
\nabla_{\dpartial_i} \left( \partial_j + f_j \partial_s + \uf_j \partial_{\us}  \right)
\\
=&
\nabla_{\dpartial_i}  \partial_j + f_j \nabla_{\dpartial_i}  \partial_s + \uf_j \nabla_{\dpartial_i}  \partial_{\us} + f_{ij} \partial_s + \uf_{ij} \partial_{\us}
\nonumber
\\
=&
f_{ij} \partial_{s}  + \uf_{ij} \partial_{\us} + \nabla_{\partial_i} \partial_j
+
\left( f_i \nabla_{\partial_j} \partial_s + f_j \nabla_{\partial_i} \partial_s \right)
+
\left( \uf_i \nabla_{\partial_j} \partial_{\us} + \uf_j \nabla_{\partial_i} \partial_{\us} \right)
\nonumber
\\
&
+
f_i f_i \nabla_{\partial_s} \partial_s + \uf_i \uf_j \nabla_{\partial_{\us}} \partial_{\us}
+
(f_i \uf_j + f_j \uf_i ) \nabla_{\partial_s} \partial_{\us}. 
\nonumber
\end{align}
Then the lemma follows from lemma \ref{lem A.1}.
\end{proof}

\begin{lemma}\label{lem A.3}
Let $\Sigma$ be a spacelike surface parametrised by $(f,\uf)$ in the double null coordinate system $\{s,\us,\vartheta\}$. $\{ \dpartial_1 , \dpartial_2\}$ is the tangential coordinate vector frame on $\Sigma$. $\{\dL, \duL\}$ is the null normal frame on $\Sigma$ constructed as in section \ref{section 3}. We have
\begin{align}
\nonumber
\nabla_{\dpartial_i} \duL
=&
\Pi_{i \duL}^k \partial_k + \Pi_{i\duL}^{\uL} \cdot \uL + \Pi_{i\duL}^L \cdot L
\end{align}
where
\begin{align}
\nonumber
\Pi_{i\duL}^k
=&
\uchi_i^k + \uvareps \chi_i^k + f_i \left[ -\uchi(b)^k + \uchi(\vec{\vareps}^k -2\Omega^2 \uvareps \eta^k - \uvareps \chi(b)^k - \slashnabla_{\vec{\uvareps}} b^k   \right]
\\
&
+
\uf_i \left[ \chi(\uvareps)^k - 2\Omega^2 \ueta^k \right] + \slashnabla_{\dpartial_i} \vec{\uvareps}^k,
\nonumber
\\
\nonumber
\Pi_{i\duL}^{\uL}
=&
\eta_i - \frac{1}{2} \Omega^{-2} \chi(\uvareps)_i + f_i \left[ 2\uomega - \eta(b) + \eta(\vec{\uvareps}) + \frac{1}{2} \Omega^{-2} \chi(b,\vec{\uvareps}) \right],
\\
\nonumber
\Pi_{i\duL}^L
=&
\uvareps \ueta_i - \frac{1}{2} \Omega^{-2} \uchi(\vec{\uvareps})_i + f_i \left[ -\uvareps \ueta(b) + \frac{1}{2} \Omega^{-2} \uchi(b,\uvareps)  \right]
+ \uf_i \left[ 2\omega \uvareps + \ueta( \uvareps)   \right] + \dpartial_i \uvareps.
\end{align}
\end{lemma}
\begin{proof}
\begin{align}
\nonumber
\nabla_{\dpartial_i} \duL
=&
\nabla_{\dpartial_i} \uL + \uvareps \nabla_{\dpartial_i} L + \uvareps^k \nabla_{\dpartial_i} \partial_k + \dpartial_i \uvareps \cdot L + \dpartial_i \uvareps^k \cdot \partial_k.
\end{align}
Then lemma follows from proposition \ref{SCACS} and
\begin{align}
\nonumber
\dpartial_i=\left(\delta_i^j - f_i b^j \right) \partial_j + f_i \underline{L} + \underline{f}_i L = B_i^j \partial_j + f_i \underline{L} + \underline{f}_i L.
\end{align}
\end{proof}

Now we can prove proposition \ref{pro 4.1}.
\begin{proof}[Proof of Proposition \ref{pro 4.1}]
\begin{align}
\nonumber
\duchi(\dpartial_i,\dpartial_j) 
=& 
- g\left( \nabla_{\dpartial_i} \dpartial_j, \duL  \right)
\\
=&
- g \left( \Pi_{ij}^k \partial_k +  \Pi_{ij}^L  L + \Pi_{ij}^{\uL} \uL,  \uL + \uvareps L + \uvareps^l \partial_l \right)
\nonumber
\\
=&
-\slashg_{kl} \Pi_{ij}^k \uvareps^l - 2\Omega^2 \Pi_{ij}^L - 2 \Omega^2 \Pi_{ij}^{\uL} \uvareps.
\nonumber
\end{align}
Then the formula of $\duchi_{ij}$ follows from lemma \ref{lem A.2}. Note that the formulae in lemma \ref{lem A.2} involve quantities $\slashGamma_{ij}^{\phantom{ij}k}, f_{ij}, \uf_{ij}$ which are not tensorial. Thus in order to simplify the calculations and obtain the tensorial formula of $\duchi_{ij}$, we can choose a geodesic coordinate system at a point such that $\slashGamma_{ij}^{\phantom{ij}k}$ vanishes and $f_{ij}=\slashnabla^2_{ij} f$, $\uf_{ij} = \slashnabla^2_{ij} \uf$ at the point.

Similarly for $\dchi_{ij}$, where
\begin{align}
\nonumber
\dchi(\dpartial_i, \dpartial_j)
=& 
- g\left( \nabla_{\dpartial_i} \dpartial_j, \dL  \right)
\\
=&
- g \left( \Pi_{ij}^k \partial_k +  \Pi_{ij}^L  L + \Pi_{ij}^{\uL} \uL,  L + \vareps \uL + \vareps^l \partial_l \right)
\nonumber
\\
=&
-\slashg_{kl} \Pi_{ij}^k \vareps^l - 2\Omega^2 \Pi_{ij}^{\uL} - 2 \Omega^2 \Pi_{ij}^{L} \vareps.
\nonumber
\end{align}
For the torsion $\deta_i$, we have
\begin{align}
\nonumber
\deta(\dpartial_i)
=&
\frac{1}{2} \dot{\Omega}^{-2} g \left( \nabla_{\dpartial_i} \duL, \dL  \right)
\\
=&
\frac{1}{2} \dot{\Omega}^{-2} g \left( \Pi_{i\duL}^{k} \partial_k + \Pi_{i\duL}^{\uL} \cdot \uL + \Pi_{i \duL}^L \cdot L, L + \vareps \uL + \vareps^l \partial_i \right)
\nonumber
\\
=&
\frac{1}{2} \Omega^{-2} \left( \slashg_{kl} \Pi_{i\duL}^k \vareps^l + 2 \Omega^2 \Pi_{i \duL}^{\uL} + 2 \Omega^2 \vareps \Pi_{i\duL}^L \right).
\nonumber
\end{align}
Then the formula of $\deta$ follows from lemma \ref{lem A.3}.
\end{proof}

\end{document}